\documentclass[11pt]{amsart}

\theoremstyle{plain}
\newtheorem{theorem}{Theorem}[section]
\newtheorem{lemma}[theorem]{Lemma}
\newtheorem{proposition}[theorem]{Proposition}
\newtheorem{corollaire}[theorem]{Corollary}

\def\cal{\mathcal}
\theoremstyle{definition}

\theoremstyle{remark}

\theoremstyle{definition}
\newtheorem{definition}{Definition}[section]

\theoremstyle{remark}
\newtheorem{remark}{Remark}

\newcommand{\M}{{\cal M}}
\newtheorem{example}{Example}

\begin{document}

\def\C{{\Bbb C}}
\def\N{{\Bbb N}}
\def\K{{\Bbb K}}
\newfont{\complg}{msbm10 scaled\magstep2}
\newcommand{\Cg}{\mbox{\complg C}}
\newcommand{\di} {{\ \rm div}}
\newcommand{\Ima} {{\Im}{\rm m}\,}
\newcommand{\Ra} {{\Re}{\rm e}\,}
\bibliographystyle{plain}
\newcommand{\z}{\zeta}
\title{Formal biholomorphic maps of real analytic
hypersurfaces}
\author{{\sc Nordine Mir}}
\address{Universit\'e de  Rouen, Laboratoire de Math\'ematiques
Rapha\"el Salem, CNRS, Site Colbert  76821 Mont Saint Aignan
France }
\email{Nordine.Mir@univ-rouen.fr}
\commby{Rouen}
\date{}


\bibliographystyle{plain}

\begin{abstract} Let $f :  (M,p) \rightarrow
(M',p')$ be a formal biholomorphic mapping between two
germs of real analytic hypersurfaces in $\C^n$,
$p'=f(p)$.  Assuming the source manifold to be minimal at $p$, we prove the convergence of the
so-called reflection function associated to $f$.  As a consequence, we derive the convergence of formal biholomorphisms
 between real analytic minimal
holomorphically nondegenerate hypersurfaces. Related results on
partial convergence of formal biholomorphisms are
also obtained. 
\end{abstract}

\maketitle

\section{Introduction} \setcounter{equation}{0} 
A formal (holomorphic) mapping $f : (\C^n,p) \rightarrow
(\C^n,p')$, $p,p'\in \C^n$, $n\geq 1$, is a vector
$(f_1(z),\ldots,f_n(z))$ where each $f_j(z)\in \C[[z-p]]$, the
ring of formal holomorphic power series in $z-p$, and
$f(p)=p'$. The mapping $f$ is called a {\it formal biholomorphism}
if its Jacobian does not vanish at $p$. If
$M,M'$ are two smooth real real-analytic hypersurfaces in $\C^n$
through
$p$ and
$p'$ respectively, we say that a formal mapping $f$ as above
sends $M$ into $M'$ if 
$\rho'(f(z),\overline{f(z)})=a(z,\bar{z})\rho (z,\bar{z})$, where
$\rho, \rho'$ are local real-analytic defining functions for
$(M,p)$ and $(M',p')$ respectively and $a\in
\C[[z-p,\bar{z}-\bar{p}]]$. In this paper we study the convergence (and partial convergence) of formal biholomorphic
mappings between germs of real analytic hypersurfaces in $\C^n$ in terms of
optimal and natural geometric
conditions on the source and target manifolds.

A natural geometric
condition which appears in this regularity problem is the concept
of {\it holomorphic nondegeneracy}.  Following Stanton,
a real analytic hypersurface
$M\subset \C^n$ is called holomorphically nondegenerate if, near
any point in $M$, there is no non-trivial holomorphic vector
field, with holomorphic
coefficients, tangent to $M$ near that point
\cite{St, BERb}.  Baouendi and Rothschild \cite{BR}
recognized the importance of such a condition and used it to characterize 
those real algebraic hypersurfaces for which any biholomorphic
self-map must be algebraic (see also
\cite{BHR}).  In this paper we  establish a
 similar statement for formal biholomorphic mappings of real
analytic hypersurfaces (Theorem \ref{cor2} below), namely that
any formal  biholomorphism between germs of  real
analytic minimal and holomorphically
nondegenerate hypersurfaces is convergent. This statement will in
fact be a consequence of our main result, Theorem \ref{th1}, where
a description of the analyticity properties of formal
biholomorphic maps of minimal real analytic
hypersurfaces is given. An application of this theorem to partial
convergence of such maps is also given in \S \ref{secpc}.

  The study of the convergence of formal mappings
between real analytic CR submanifolds goes back to
Chern and Moser \cite{CM}, who established the
convergence of formal biholomorphisms between real analytic 
Levi-nondegenerate hypersurfaces. More recently, Baouendi,
Ebenfelt and Rothschild \cite{BER2, BER3, BER4} addressed this
problem in more general situations. In particular, the
following two facts follow from their work:

 i) Given a generic
minimal real analytic connected holomorphically nondegenerate
submanifold
$M$ in
$\C^N$, 
$N\geq 2$, there exists a proper real analytic subvariety
$S\subset M$, such that for any point $p\in M\setminus S$, any
formal biholomorphism sending $(M,p)$ onto another germ of a
generic real analytic submanifold in $\C^N$, must be
convergent.

 ii) If $M$ is a connected generic holomorphically
degenerate submanifold in $\C^N$, then, for any point $p\in M$,
there exists a formal biholomorphism sending $(M,p)$ into
itself which does not converge.

 In view of these facts and other related results in
the mapping problems (\cite{BERb}), it appears likely that the
 subvariety $S$ in i) could be taken to be empty. In this
paper, we actually prove such a result in the case of
hypersurfaces.  Our proof is based on an analysis of a so-called reflection
function associated to the mapping and the hypersurfaces, which
was used in many other situations ({\it cf.} \cite{H, M2}).
We should also mention that, recently in \cite{M3}, we showed
that the real analytic subvariety $S$ in i) can also be taken
to be empty when the generic submanifold $M$ is assumed to be
real algebraic (and with no such assumption on the target
manifold).
\section{Statement of main results}\label{sec1}
\setcounter{equation}{0}

Let $(M',p')\subset \C^n$, $n\geq 2$, be a germ
at $p'$ of a smooth real real-analytic hypersurface.
Let $\rho'=\rho' (\zeta,\bar{\zeta})$ be a real
analytic defining function for
$M'$ near $p'$, i.e.
$$M'=\{\zeta \in (\C^n,p'):\rho' (\zeta,\bar{\zeta)}=0\}.$$ After
complexification of $\rho'$, one defines the so-called invariant
{\it Segre varieties} attached to $M'$ by
$$Q'_{\omega}=\{\zeta \in (\C^n,p'):\rho'(\zeta,\bar{\omega})=0
\},$$ for $\omega$ close to $p'$. We can assume, without loss of
generality, that $p'=0$ and $\frac{\partial \rho'}{\partial
\zeta_n}(0)\not =0$, $\zeta=(\zeta',\zeta_n)\in \C^{n-1}\times
\C$. Thus, the implicit function theorem exhibits any Segre
variety as a graph of the form
$$Q'_{\omega}=\{\zeta \in
(\C^n,0):\zeta_n=\Phi'(\bar{\omega},\zeta')\},$$ where $\Phi'$ is
a holomorphic function in its arguments in a neighborhood of
$0\in \C^{2n-1}$ and such that $\Phi'(0)=0$. Equivalently, this Segre
variety can be defined by
$$Q'_{\omega}=\{\zeta \in
(\C^n,0):\bar{\zeta}_n=\bar{\Phi}'(\omega,\bar{\zeta}')\}.$$
Here, we have used the following notation. If $g=g(x)$ is some formal
holomorphic power series in $\C[[x]]$,
$x=(x_1,\ldots,x_k)$, 
$\bar{g}$ is the formal holomorphic power series
obtained by taking the complex conjugates of the coefficients of
$g$. (This convention of notation will be used throughout the
paper.) Our main result is the following.

\begin{theorem}\label{th1} Let $f : (M,0)\rightarrow
 (M',0)$ be a formal biholomorphism
between two germs at 0 of  real-analytic
hypersurfaces in
$\C^n$.  Assume, furthermore, that
$M$ is minimal at $0$. 
Then, the formal holomorphic map
$$\C^n\times \C^{n-1}\ni (z,\lambda)\mapsto \bar{\Phi}'(f(z),\lambda)$$
is convergent.  \end{theorem} 
Such a result has several applications.  One of its main applications
lies in the following  theorem mentioned in the introduction. 

\begin{theorem}\label{cor2} Any formal biholomorphic mapping
between germs of minimal, holomorphically
nondegenerate, real analytic hypersurfaces in $\C^n$ is
convergent. 
\end{theorem}

 As explained in the
introduction, the interest in such a result lies in the fact
that, in view of \cite{BER2}, the condition of holomorphic
nondegeneracy is optimal for the class of formal biholomorphisms. 
However, it remains an open problem to decide whether or not the
condition of minimality is necessary in Theorem \ref{cor2}. Another
application of Theorem \ref{th1} deals with partial convergence of
formal biholomorphisms. For this, we refer the reader to \S
\ref{secpc}.\begin{remark}  After this work was
completed, I received a preprint by J.  Merker,
``Convergence of formal biholomorphic mappings between
minimal holomorphically nondegenerate real analytic
hypersurfaces'', in which a similar statement to
Theorem
\ref{cor2} is given.  In that preprint, Theorem
\ref{th1} above is also stated for the special case where $M'$ is
a rigid and polynomial hypersurface.\end{remark}   
\section{Definitions and notations}\label{sec2}
\setcounter{equation}{0} \subsection{Real analytic
hypersurfaces.}\label{ssec21} Let $M$ be a smooth real
real-analytic hypersurface in $\C^n$, $n\geq 2$.  Since the
situation is purely local, we shall always work near a point
$p\in M$, which will be assumed, without loss of generality, to
be the origin.  Let $\rho$ be a real analytic defining function
for $M$ near 0 i.e.  $$M=\{z\in (\C^n,0):\rho (z,\bar{z})=0\},$$
with $d\rho
\not = 0$ on $M$.  The complexification $\M$ of $M$ is the
complex hypersurface through 0 in $\C^{2n}$ given by
$$\M=\{(z,w)\in (\C^{2n},0):\rho (z,w)= 0 \}.$$ We shall assume,
without loss of generality, that the coordinates $z=(z',z_n)\in
\C^n$ are chosen so that $\displaystyle
\frac{\partial
\rho}{\partial z_n}(0)\not =0$.  In this case, we define the
following holomorphic vector fields tangent to
$\M$ 
\begin{equation}\label{eqvf}
{\cal L}_j=\frac{\partial \rho}{\partial
w_n}(z,w)
\frac{\partial}{\partial w_j}-\frac{\partial
\rho}{\partial w_j}(z,w)\frac{\partial}{\partial w_n},\
j=1,\ldots,n-1,
\end{equation} which are the complexifications of the usual
(0,1) vector fields tangent to $M$.  We recall that for a point
$w$ near 0, its associated Segre surface is the complex
hypersurface defined by $Q_{w}=\{z\in (\C^{n},0):\rho (z,
\bar{w})=0\}$.  Observe that by the complex analytic implicit
function theorem, each Segre variety can be described as a graph
of the form
$$Q_{w}=\{z\in (\C^{n},0):z_n=\Phi (\bar{w},z')\},$$ $\Phi$
denoting a convergent power series in some neighborhood of the
origin in $\C^{2n-1}$ satisfying the relations $\Phi (0)=0$ and
\begin{equation}\label{eqm1} \Phi (w',\bar{\Phi}(z,w'),z')\equiv
z_n,\quad (z,w')\in \C^n\times \C^{n-1}.
\end{equation} Equation (\ref{eqm1}) is a consequence of the fact
that $M$ is a {\it real} hypersurface.  The coordinates $z$ are
said to be {\it normal} with respect to $M$ if the additional
condition
$$\Phi(w',w_n,0)=\Phi (0,w_n,z')\equiv w_n,\quad (w',w_n,z')\in
\C^{n-1}\times \C \times \C^{n-1},$$ holds. It is well-known that
given a real-analytic hypersurface through the
origin, one can always construct such coordinates \cite{CM}. 
Thus, from now on and for simplicity, we will always assume that the
$z$-coordinates are chosen to be normal for
the manifold $M$.

 The real analytic hypersurface $M$ is called {\it minimal} at 0 
 (in the sense of Tr\'epreau and Tumanov), or, equivalently of {\it
finite type} (in the sense of Kohn and Bloom-Graham) if it does not
contain any complex-analytic hypersurface through 0.  To use such a
nondegeneracy condition, we will need the {\it Segre set mappings}
associated to $M$ (see
\cite{BERb}) up to order 3, which, in normal coordinates, are the
following three maps : 
\begin{equation}\label{market}
\begin{aligned}(\C^{n-1},0)\ni
z'\mapsto v_1(z')&=(z',0)\in \C^n,\\
 (\C^{2n-2},0)\ni
(z',\xi)\mapsto v_2(z',\xi)&=(z',\Phi (\xi,0,z'))\in \C^n,\\ 
(\C^{3n-3},0)\ni (z',\xi,\eta)\mapsto v_3(z',\xi,\eta)&=(z',\Phi
(\xi,\bar{\Phi}(\eta,0,\xi),z'))\in \C^n.
\end{aligned}  
\end{equation}These maps are of
fundamental importance since they parametrize the so-called {\it
Segre sets} (up to order three) associated to $M$.  Moreover, the
interest in such Segre sets or Segre set mappings lies in the fact
that the minimality assumption is equivalent to the fact that the
generic rank of $v_2$ (and also $v_3$) equals $n$ (see \cite{BERb}). 
This will be useful for the proof of Theorem \ref{th1}.\par
\smallskip All the notations introduced in this section will
be used for the source hypersurface $M$.  For the target real
analytic hypersurface $M'$, we shall use the notations introduced in
\S \ref{sec1}, before the statement of Theorem \ref{th1}. 
In particular, we denote
by $z$ the coordinates in the source space and by $\zeta$ the
coordinates at the target space.

 \subsection{Some commutative algebra.} We recall here, for the
reader's convenience, some basic definitions,  needed in \S
\ref{secpc}, about regular local rings and their ideals. All
these definitions can be found for instance in \cite{L, T}.

Let
$A$ be a Noetherian ring. If $I$ and $J$ are two ideals of $A$,
we use the notation $I<J$ to mean that $I\subset J$ and $I\not
=J$. Given a prime ideal $I\subseteq A$, the {\it height} of $I$
is defined by the formula
$$h(I)={\rm max}\{k\in \N:\{0\}<I_1<\ldots<I_k=I\},$$ where
$I_{1},\ldots,I_k$ are prime ideals of $A$. If
$J$ is any ideal of $A$, we define the height of $J$ by the formula
$$h(J)={\rm inf}\{h(I):J\subset I,\ I\  {\rm prime\ ideal\ of}\
A\}.$$ If $A$ is furthermore assumed to be a local ring, one
defines the {\it Krull dimension} of $A$ to be the height of its
maximal ideal. Observe that if $A$ is a Noetherian local ring and
if $I$ is a proper ideal of $A$, the quotient ring $A/I$ is also
a Noetherian local ring. This allows one to consider the
Krull dimension of such a ring. Finally, a Noetherian local
ring is said to be {\it regular} if its maximal ideal has
$\delta$ generators, where $\delta$ is the Krull dimension of the
ring $A$. The rings of formal holomorphic power series or
convergent power series in $p$ indeterminates, $p\in \N^*$,
are regular rings of Krull dimension $p$ \cite{T}. 
\section{Two convergence results}\label{sec-1}
\setcounter{equation}{0}
In this section, we 
first state and prove a convenient lemma which will be used twice in
the paper. This lemma may already be known.

\begin{lemma}\label{lem00} Let $(u_i(t))_{i\in I}$ be a family of
convergent power series in $\C\{t\}$, 
$t=(t_1,\ldots,t_q)$, $q\in \N^*$. Let also $({\cal
K}_{i}(\varsigma))_{i\in I}$ be a family of convergent power
series in $\C\{\varsigma\}$, $\varsigma
=(\varsigma_1,\ldots,\varsigma_r)$, $r\in \N^*$. Assume
that:
\begin{enumerate}
\item [(i)] There exists $R>0$ such that the radius of convergence
of any ${\cal K}_i$, $i\in I$, is at least $R$. 
\item [(ii)] For all
$\varsigma \in \C^r$ with $|\varsigma|<R$, $|{\cal
K}_i(\varsigma)|\leq C_i$, with $C_i>0$.
\item [(iii)] There exists
$V(t)=(V_1(t),\ldots,V_r(t))\in (\C[[t]]
)^r$, $V(0)=0$, such that
$({\cal K}_i\circ V)(t)=u_i(t)$ (in $\C[[t]]$) for all $i\in I$.
\end{enumerate}
Then, there exists $R'>0$ such that the radius of convergence of
any $u_i$, $i\in I$, is at least $R'$ and such that for all $t\in
\C^q$ with $|t|<R'$, $|u_i(t)|\leq C_i$. 
\end{lemma}
\begin{proof}  By Artin's approximation
theorem
\cite{A}, one can find a convergent power series
mapping 
$\vartheta (t)=(\vartheta_1(t),\ldots,\vartheta_r(t))\in
(\C\{t\})^r$ such that $\vartheta (0)=0$ and for all 
$i\in I$, $({\cal K}_i \circ \vartheta) (t)=u_i(t)$ in $\C\{t\}$.
Let $R'>0$ so that if $|t|<R'$ then 
$|\vartheta (t)|<R$. Since the radius of convergence of the
${\cal K}_i$, $i\in I$, is at least $R$, the radius  of
convergence of any ${\cal K}_i\circ \vartheta$, $i\in I$, is at
least $R'$. Thus, the family $(u_i(t))_{i\in I}$ has a radius of
convergence at least equal to $R'$ and for all $t\in \C^q$ with
$|t|<R'$, one has $|u_i(t)|= |({\cal K}_i\circ \vartheta) (t)|\leq
C_i$, $i\in I$.\end{proof}
 To derive Theorem \ref{cor2} from Theorem \ref{th1},
we will need the following consequence of Artin's approximation
theorem, which is contained for instance in \cite{N}.  (See also
\S \ref{secpc} for another formulation of such a result.) 
\begin{proposition}\label{prop22} Let
$R(x,y)=(R_1(x,y),\ldots,R_r(x,y))\in (\C\{x,y\})^r$, $x\in
\C^q$, $y\in \C^r$, $q,r\in \N^*$.  Let
$g(x)=(g_1(x),\ldots,g_r(x))\in (\C[[x]])^r$ satisfy  $R(x,g(x))=0$.  If ${\rm
det}\big(\displaystyle \frac{\partial R}{\partial y}(x,g(x))\big)\not
\equiv 0$ in $\C[[x]]$, then $g(x)$ is convergent. 
\end{proposition} \begin{proof}  We reproduce here
the arguments of \cite{N}.  Write
\begin{equation}\label{eq56} R(x,y)-R(x,z)=Q(x,y,z)\cdot(y-z)
\end{equation} where $Q$ is an $r\times r$ complex-analytic matrix
such that
$Q(x,y,y)=\displaystyle \frac{\partial R}{\partial y}(x,y)$; i.e.\
 $Q(x,y,z)=\displaystyle
\int_0^1\displaystyle \frac{\partial R}{\partial
y}(x,ty+(1-t)z)dt$. By assumption, we know that we 
have
${\rm det\,} Q(x,g(x),g(x))\not \equiv 0$.  This
implies that one can find an integer $k_{g}$ such
that if $H(x)$ is any formal power series which
agrees up to order $k_{g}$ with $g$ then ${\rm
det\,}Q(x,g(x),H(x))\not \equiv 0$. For this
integer $k_{g}$, according to  Artin's
approximation theorem, one can find a convergent
power series
$H_0(x)$ satisfying $R(x,H_0(x))=0$ and agreeing with $g(x)$ up to
order
$k_{g}$.  By (\ref{eq56}), we get
$Q(x,g(x),H_0(x))\cdot(g(x)-H_0(x))\equiv 0$ in $\C[[x]]$.  Since
${\rm det\,}Q(x,g(x),H_0(x))\not \equiv 0$,
we obtain $g(x)=H_0(x)$ and thus $g$ is
convergent.\end{proof} 
\section{The reflection principle}\label{sec3}
\setcounter{equation}{0} Let $f:(M,0)\rightarrow (M',0)$ be a
formal biholomorphic mapping between germs at 0 of real analytic hypersurfaces
in $\C^n$. We shall use the notations introduced in
\S \ref{sec1} and \S\ref{sec2}. In particular, we
denote $f=f(z)=(f_1(z),\ldots,f_n(z))=(f'(z),f_n(z))$ in the
$\zeta$-coordinates. In this section, we shall make no further
assumptions on $M$ and $M'$.

As in \cite{M2, M3}, we define
the following formal holomorphic power series 
\begin{equation}\label{eqR} \C^n\times \C^{n-1} \ni
(z,\lambda)\mapsto {\cal
R}(z,\lambda):=\bar{\Phi}'(f(z),\lambda).
\end{equation} The goal of this section
is to prove the following proposition.
\begin{proposition}\label{prop1} Let $f:(M,0)\rightarrow (M',0)$ be a formal biholomorphism
between germs at 0 of real analytic hypersurfaces in $\C^n$, and
${\cal R}$ defined by {\rm (\ref{eqR})}. Then for any multi-index
$\gamma
\in
\N^{n}$, the formal holomorphic map $$(\C^{2n-2},0)\ni
(z',\lambda) \mapsto \left(\partial_{z}^{\gamma}{\cal
R}(z,\lambda)\right){\big|_{z=v_1(z')}}$$ is convergent in some
neighborhood
$V_{\gamma}$ of $0\in \C^{2n-2}$. Here, $v_1$ is the first Segre set
mapping for $M$ as defined in {\rm (\ref{market})}. 
\end{proposition} 
Before proceeding to the proof of Proposition \ref{prop1}, we need a
preliminary lemma (Lemma
\ref{lem1} below). Since $f$ maps formally $M$ into $M'$, there exists
$a(z,\bar{z})\in
\C[[z,\bar{z}]]$ such that
$$\overline{f_n(z)}-\bar{\Phi}'(f(z),\overline{f'(z)})=a(z,\bar{z})
\rho(z,\bar{z}),\ {\rm in}\ \C[[z,\bar{z}]].$$ Equivalently, we
have
\begin{equation}\label{eq4000}
\bar{f}_n(w)-\bar{\Phi}'(f(z),\bar f'(w))=a(z,w) \rho (z,w),\
{\rm in}\
\C[[z,w]].
\end{equation}
We write $\bar\Phi'_{\lambda^{\alpha}}(\omega,\lambda)$ for
$\partial_{\lambda}^{\alpha}\bar\Phi'(\omega,\lambda)$. By
applying the vector fields
${\cal L}_j$,
$j=1,\ldots,n-1$, as defined by (\ref{eqvf}), to (\ref{eq4000})
and using the fact that $f$ is invertible, one obtains the
following known statement (see
\cite{BER3} for instance). 

\begin{lemma}\label{lem1} Under the assumptions of Proposition
$\ref{prop1}$, one has, for any multindex
$\alpha
\in
\N^{n-1}$, the formal power series identity $$
\bar\Phi'_{\lambda^{\alpha}}(f(z),\bar
f'(w))=\chi_{\alpha}\left((\partial^{\beta}\bar
f(w))_{|\beta|\leq |\alpha|},z,w\right),\ (z,w)\in {\cal M},$$
where each
$\chi_{\alpha}$ is a convergent power series of its arguments. 
\end{lemma} \begin{proof}[Proof of Proposition
$\ref{prop1}$.]  We write the expansion
\begin{equation}\label{eqts}
\bar\Phi'(\omega,\lambda)=\sum_{\alpha \in
\N^{n-1}}\phi'_{\alpha}(\omega)
\lambda^{\alpha}. 
\end{equation}
For the
sake of clarity, we shall first give the proof of the Proposition in
the case
$\gamma =0$.

{\sc The case $\gamma =0$.} We restrict all the
identities given by Lemma \ref{lem1} to the
$(n-1)$-dimensional subspace $$\{(0,v_1(z')):z'\in
(\C^{n-1},0)\}\subset {\cal M}.$$ This gives, for any
multiindex $\alpha \in \N^{n-1}$, 

\begin{equation}\label{eq1} \alpha!\ (\phi'_{\alpha}\circ
f\circ
v_1)(z')=\chi_{\alpha}((\partial^{\beta}\bar{f}(0))_{|\beta|\leq
|\alpha|},v_1(z'),0):=u_{\alpha}(z'),
\end{equation} 
where $\phi'_{\alpha}$ is given by (\ref{eqts}). Observe that for each
multiindex $\alpha$, $u_{\alpha}(z')$ is convergent. 

To show that
$\C^{n-1}\times \C^{n-1}\ni (z',\lambda)
\mapsto {\cal R}(v_1(z'),\lambda)$ is
convergent, we claim that it suffices to show that there exists $a>0$
and $R_0>0$ such that the radius of convergence of each $u_{\alpha}$ is
at least $a$, and such that the following Cauchy estimates hold:
\begin{equation}\label{cauchy}
\forall \alpha \in \N^{n-1},\ \forall
z'\in
\C^{n-1},\ |z'|<a,\ |u_{\alpha}(z')|\leq \alpha !\
R_0^{|\alpha|+1}.
\end{equation}  Indeed, if (\ref{cauchy}) holds, then the formal
holomorphic power series
$$ \C^{n-1}\times \C^{n-1}\ni (z',\lambda)\mapsto
{\cal R}_0(z',\lambda):=\sum_{|\alpha|=0}^{\infty}\displaystyle
\frac{u_{\alpha}(z')}{\alpha!}\lambda^{\alpha}$$ defines
a convergent power series in $B_{n-1}(0,a)\times
B_{n-1}(0,{1}/{2R_0})$. (Here and in what follows, for any $c>0$
and for any $k\in \N^*$, $B_k(0,c)$  denotes the
euclidean ball centered at 0 in $\C^k$ of radius $c$.) Moreover, by
(\ref{eq1}) and (\ref{eqR}), we have for
any multiindex
$\alpha \in \N^{n-1}$, \begin{equation}\label{eq4001}
\left[\frac{\partial^{|\alpha|} {\cal R}_0}{\partial
\lambda^{\alpha}}(z',\lambda)\right]_{\lambda=0}=\left[
\frac{\partial^{|\alpha|} {\cal R}}{\partial
\lambda^{\alpha}}(v_1(z'),\lambda)\right]_{\lambda=0} \end{equation}
in  $\C[[z']]$ and hence 
$${\cal R}_0(z',\lambda)= {\cal R}(v_1(z'),\lambda).$$ This proves
that under the assumption (\ref{cauchy}), ${\cal
R}(v_1(z'),\lambda)$ is a convergent power series in
$(z',\lambda)$. It remains to prove (\ref{cauchy}). Since
$\bar\Phi'$ is holomorphic in a neighborhood of $0\in \C^{2n-1}$,
in view of (\ref{eqts}), one can find
$\delta >0$ and a constant
$R>0$ such that for any multiindices
$\alpha \in \N^{n-1}$, $\nu \in \N^{n}$,
\begin{equation}\label{eq2} \forall \omega\in \C^n,\
|\omega|<\delta,\ 
\left|\frac{\partial^{|\nu|} \phi'_{\alpha}(\omega)}{\partial
\omega^{\nu}}\right|\leq \nu!\ R^{|\alpha|+|\nu|+1}. 
\end{equation} In view of (\ref{eq1}) and (\ref{eq2}) (in the
case $\nu=0$), we can apply Lemma \ref{lem00} to conclude that
there exists $a>0$ such that the family $(u_{\alpha}(z'))_{\alpha
\in \N^{n-1}}$ is convergent in $B_{n-1}(0,a)$ and such that
(\ref{cauchy}) holds with $R_0=R$. This finishes the proof of
Proposition
\ref{prop1} in the case
$\gamma=0$.  

{\sc The case $|\gamma|>0$.}  We proceed now to the proof
of Proposition \ref{prop1} for general $\gamma \in \N^n$.  For this, we need the following lemma.

\begin{lemma}\label{lem2} Under the assumptions of Proposition
$\ref{prop1}$, for any multiindices
$\alpha\in
\N^{n-1}$, $\gamma \in \N^{n}$, the formal holomorphic power
series $$\C^{n-1}\ni z'\mapsto \partial^{\gamma}_{z} \left(
(\phi'_{\alpha}\circ
f)(z)\right)\big|_{z=v_1(z')}$$ is convergent.  \end{lemma}
\begin{proof}[Proof of Lemma $\ref{lem2}$.]  We
prove the Lemma by induction on
$|\gamma|$ (for any multiindex $\alpha \in
\N^{n-1}$).  For $\gamma=0$, the statement follows from (\ref{eq1}), as we previously noticed.  Let
$\gamma\in \N^n$.  For $\alpha \in \N^{n-1}$, using Lemma
\ref{lem1}, we obtain for
$(z',z_n,0,z_n)\in (\M,0)$,
$$\bar\Phi'_{\lambda^{\alpha}}(f(z),\bar{f}'(0,z_n))=
\chi_{\alpha}\left((\partial^{\beta}\bar{f}(0,z_n))_{|\beta|\leq
|\alpha|},z,0,z_n\right).$$ If we apply $\partial^{\gamma}_{z}$
to this equation, we obtain
\begin{equation}\label{eq4002}
\frac{\partial^{|\gamma|}}{\partial z^{\gamma}}\left[{\cal
R}_{\lambda^{\alpha}}(z,\bar{f}'(0,z_n))\right]=
\frac{\partial^{|\gamma|}}{\partial
z^{\gamma}}\left[
\chi_{\alpha}((\partial^{\beta}\bar{f}(0,z_n))_{|\beta|\leq
|\alpha|},z,0,z_n)\right].  \end{equation} One can easily check
that this implies that there exist a polynomial ${\cal
S}_{\gamma}$ such that the left-hand side of
(\ref{eq4002}) is equal to 
\begin{multline}\label{eq4003}
{\cal R}_{z^{\gamma}\lambda^{\alpha}}
(z,\bar{f}'(0,z_n))+\\ {\cal
S}_{\gamma}\left[ \left({\cal
R}_{z^{\nu}\lambda^{\beta}}(z,\bar{f}'(0,z_n))\right)_
{|\nu|<|\gamma| \atop |\beta|\leq
|\alpha|+|\gamma|},\left(\partial^{\mu}\bar{f}(0,z_n)\right)
_{|\mu|\leq
|\gamma|}\right],  \end{multline}
where $\mu, \nu \in \N^n$, $\beta \in \N^{n-1}$. Furthermore, we
observe that the right-hand side of (\ref{eq4002}) can be
written in the form
\begin{equation}\label{eq4004}
\chi^1_{\alpha,\gamma}\left((\partial^{\beta}\bar{f}(0,z_n))_{|\beta|\leq
|\alpha|+|\gamma|},z\right), \end{equation} where
$\chi^1_{\alpha,\gamma}$ is a convergent power series.  Restricting (\ref{eq4002}), (\ref{eq4003}) and
(\ref{eq4004}) to $z=v_1(z')$, one obtains
\begin{multline}
\alpha!\ \partial^{\gamma}_{z} \left(
(\phi'_{\alpha}\circ
f)(z)\right)\big|_{z=v_1(z')}+\\
{\cal S}_{\gamma}\left[\big(\beta!\
\partial^{\nu}_{z} \left(
(\phi'_{\beta}\circ
f)(z)\right)\big|_{z=v_1(z')}\big)_{|\nu|<|\gamma| \atop 
|\beta|\leq
|\alpha|+|\gamma|},(\partial^{\mu}\bar{f}(0))_{|\mu|\leq
|\gamma|}\right]=\\
\chi^1_{\alpha,\gamma}\left((\partial^{\beta}\bar{f}(0))_{|\beta|\leq
|\alpha|+|\gamma|},z',0\right).
\end{multline} The induction hypothesis tells us
that for any multiindex $\beta \in \N^{n-1}$ and for any
multiindex $\nu \in \N^n$ such that $|\nu|<|\gamma|$, the formal
holomorphic power series 
$$z'\in \C^{n-1} \mapsto
\partial^{\nu}_{z} \left(
(\phi'_{\beta}\circ
f)(z)\right)\big|_{z=v_1(z')}$$ is
convergent.  Thus, we obtain the desired similar statement for
$\C^{n-1}\ni z' \mapsto \partial^{\gamma}_{z} \left(
(\phi'_{\alpha}\circ
f)(z)\right)\big|_{z=v_1(z')}$, for any multiindex $\alpha \in
\N^{n-1}$.\end{proof} 

We come back to the proof of
Proposition \ref{prop1}.  For all multiindices $\gamma \in
\N^{n},\ \alpha \in \N^{n-1}$, we put
\begin{equation}\label{eq3} \Psi_{\alpha,\gamma}(z'):=\alpha!\
\partial^{\gamma}_{z} \left( (\phi'_{\alpha}\circ
f)(z)\right)\big|_{z=v_1(z')}={\cal
R}_{z^{\gamma}\lambda^{\alpha}}(v_1(z'),0).
\end{equation} 
By Lemma \ref{lem2}, the $\Psi_{\alpha,\gamma}(z')$ are
convergent. We now fix $\gamma$,
$|\gamma|\geq 1$.  We want to prove that
${\cal R}_{z^{\gamma}}(v_1(z'),\lambda)$ is convergent in
some neighborhood
$V_{\gamma}$ of $0\in
\C^{2n-2}$.  For this, as in the
case $\gamma=0$, it suffices to prove that one can find
$q_{\gamma}>0$ and $R_{\gamma}>0$ such that the radius of
convergence of the family $(\Psi_{\alpha,\gamma})_{\alpha \in
\N^{n-1}}$ is at least $q_{\gamma}$ and such that the following
estimates hold:
\begin{equation}\label{eqest}
\forall \alpha \in \N^{n-1},\ \forall z'\in
\C^{n-1},\ |z'|<q_{\gamma},\ |\Psi_{\alpha,\gamma}(z')|\leq
\alpha!\ R_{\gamma}^{|\alpha|+1}.
\end{equation} 
We first observe that
for any multiindex $\nu \in \N^n$ with $|\nu|\leq |\gamma|$,
there exists a universal polynomial ${\cal P}_{\nu,\gamma}$,
such that
\begin{multline}\label{eq4006} {\cal
R}_{z^{\gamma}\lambda^{\alpha}}(z,\lambda)=
\partial^{\gamma}_{z}\left(\bar\Phi'_{\lambda^{\alpha}}
(f(z),\lambda)\right)=\\
\sum_{|\nu|\leq |\gamma|} {\cal
P}_{\nu,\gamma}\left((\partial^{\beta}f(z))_{1\leq |\beta|\leq
|\gamma|}\right)\bar\Phi'_{
\omega^{\nu}\lambda^{\alpha}}(f(z),\lambda).  \end{multline} This
means in particular that the polynomials ${\cal P}_{\nu,\gamma}$,
$|\nu|\leq |\gamma|$, are independent of $\alpha$. Putting
$\lambda=0$ and $z=v_1(z')$ in (\ref{eq4006}),  we obtain
\begin{multline}\label{eq4} {\cal
R}_{z^{\gamma}\lambda^{\alpha}}(v_1(z'),0)=\\
\alpha!\
\sum_{|\nu|\leq |\gamma|} {\cal
P}_{\nu,\gamma}\left(((\partial^{\beta}f)(
v_1(z')))_{1\leq |\beta|\leq
|\gamma|}\right)\left(\frac{\partial^{|\nu|}
\phi'_{\alpha}}{\partial \omega^{\nu}}\circ f\circ
v_1\right)(z').  \end{multline} Recall that $\gamma$ is fixed. 
For $\alpha \in \N^{n-1}$, consider the convergent
power series of the variables
$((\Lambda_{\beta})_{1\leq |\beta|\leq |\gamma|},\omega)$ defined by
$$h_{\alpha,\gamma}((\Lambda_{\beta})_{1\leq |\beta|\leq
|\gamma|},\omega):= \alpha!\ \sum_{|\nu|\leq |\gamma|} {\cal
P}_{\nu,\gamma}\left((\Lambda_{\beta}+\partial^{\beta}f(0))_{1\leq
|\beta|\leq |\gamma|}\right)\frac{\partial^{|\nu|}
\phi'_{\alpha}}{\partial
\omega^{\nu}}(\omega).$$
Let $r(\gamma)=n\ {\rm Card}
\{\nu
\in \N^n:1\leq |\nu|\leq |\gamma|\}$. In view of (\ref{eq2}),
each $h_{\alpha,\gamma}$, $\alpha \in \N^{n-1}$, is
convergent in $B_{r(\gamma)}(0,1)\times B_{n-1}(0,\delta)$. 
Moreover, for
$((\Lambda_{\beta})_{1\leq |\beta|\leq |\gamma|},\omega)\in
B_{r(\gamma)}(0,1)\times B_{n-1}(0,\delta)$, we have the following estimates
$$|h_{\alpha,\gamma}((\Lambda_{\beta})_{1\leq |\beta|\leq
|\gamma|},\omega)|\leq \alpha !\ \sum_{|\nu|\leq |\gamma|} |{\cal
P}_{\nu,\gamma}\left((\Lambda_{\beta}+\partial^{\beta}f(0))_{1\leq
|\beta|\leq |\gamma|}\right)\!\!|\ \nu !\
R^{|\alpha|+|\nu|+1}.$$ Put
\begin{multline}
C_{\gamma}:={\rm sup}\{|{\cal
P}_{\nu,\gamma}\left((\Lambda_{\beta}+\partial^{\beta}f(0))_{1\leq
|\beta|\leq |\gamma|}\right)\!|:\\|\nu|\leq |\gamma|,\
(\Lambda_{\beta}^{n})_{1\leq |\beta|\leq |\gamma|}\in
B_{r(\gamma)}(0,1)\}.
\end{multline}  This implies that in
$B_{r(\gamma)}(0,1)\times B_{n-1}(0,\delta)$, the following
estimates hold for some suitable constant $C^1_{\gamma}$:
$$|h_{\alpha,\gamma}((\Lambda_{\beta})_{1\leq |\beta|\leq
|\gamma|},\omega)|\leq \alpha !\ C_{\gamma}\sum_{|\nu|\leq
|\gamma|} \nu !R^{|\alpha|+|\nu|+1}\leq C_{\gamma}^1\ \alpha !\ 
R^{|\alpha|+|\gamma|+1}.$$  From this, we see that there exists
$R_{\gamma}>0$ such that for $((\Lambda_{\beta})_{1\leq |\beta|\leq
|\gamma|},\omega)\in B_{r(\gamma)}(0,1)\times B_{n-1}(0,\delta)$,
\begin{equation}\label{eqr1}
|h_{\alpha,\gamma}((\Lambda_{\beta})_{1\leq |\beta|\leq
|\gamma|},\omega)|\leq \alpha !\ R_{\gamma}^{|\alpha|+1}.
\end{equation} In view of (\ref{eq3}) and (\ref{eq4}), we have
for any multiindex $\alpha \in \N^{n-1}$,
$$h_{\alpha,\gamma}\left(((\partial^{\beta}f)
(v_1(z'))-\partial^{\beta}f(0))_{1\leq |\beta|\leq
|\gamma|},(f\circ v_1)(z')\right)
=\Psi_{\alpha,\gamma}(z'),$$ as formal power series in $z'$.
Thus, in view of (\ref{eqr1}), we are in a position to apply
Lemma
\ref{lem00} to conclude that there exists $q_{\gamma}>0$ such
that the family $(\Psi_{\alpha,\gamma}(z')) _{\alpha \in
\N^{n-1}}$ is convergent in $B_{n-1}(0,q_{\gamma})$, in which, moreover, 
the desired estimates (\ref{eqest}) hold. This
implies that ${\cal R}_{z^{\gamma}}(v_1(z'),\lambda)\in
\C\{z',\lambda \}$. The proof of Proposition \ref{prop1} is
thus complete.\end{proof}
\begin{remark} It is clear that Proposition
\ref{prop1} still holds in higher codimension with the same
proof. More precisely, the following holds. Let
$M,M'$ be two germs through the origin in $\C^n$,
$n\geq 2$, of smooth real real-analytic generic submanifolds
of CR dimension $N$ and of real codimension $d$.  Let $f :
(M,0)\rightarrow (M',0)$ be a  formal biholomorphic map.
Assume that the coordinates at the target space
$\zeta=(\zeta',\zeta^*)\in \C^N
\times \C^d$ are chosen so that, near the origin,
$M'=\{(\zeta',\zeta^*)\in
(\C^n,0):\bar{\zeta}^{*}=\bar{\Phi}'(\zeta,\bar{\zeta}')\}$,
for some
$\C^d$-valued holomorphic map
$\bar\Phi'=(\bar\Phi'_1,\ldots,\bar\Phi'_d)$ near $0\in \C^{2N+d}$.  Assume
also that the coordinates at the source space are chosen to be
normal coordinates for $M$. Then, if we define $\C^{n}\times
\C^N\ni (z,\lambda)\mapsto {\cal
R}(z,\lambda):=\bar\Phi'(f(z),\lambda)\in \C^d$,  the
following holds. For any multiindex $\gamma \in \N^n$, the
formal holomorphic power series mapping
$$\C^N\times \C^N \ni (z',\lambda)\mapsto {\cal
R}_{z^{\gamma}}((z',0),\lambda)\in \C^d$$ is convergent.
\end{remark}

\section{Proofs of Theorem \ref{th1} and Theorem
\ref{cor2}}\label{sec4} \setcounter{equation}{0} For the proof
of Theorem \ref{th1}, we first need to prove a
lemma (also used in \cite{M3}) which will allow us to bypass
the second Segre set and to work directly on the third Segre
set. 
\begin{lemma}\label{lem3} Let ${\cal
T}(x,u)=({\cal T}_1(x,u),\ldots,{\cal T}_r(x,u))\in
(\C[[x,u]])^r$,
$x\in \C^q$, $u\in \C^s$, with ${\cal T}(0)=0$.  Assume that
${\cal T}(x,u)$ satisfies an identity in the  ring
$\C[[x,u,y]]$, $y\in \C^q$, of the form
$$\varphi ({\cal T}(x,u);x,u,y)=0,$$
where $\varphi \in
\C[[W,x,u,y]]$ with $W\in \C^r$.  Assume, furthermore,
that for any multi-index $\beta \in \N^q$, the formal power
series $\displaystyle
\left[\frac{\partial^{|\beta|}\varphi}{\partial y^{\beta}}
(W;x,u,y)\right]_{y=x}$ is convergent,
i.e.  belongs to $\C\{W,x,u\}$.  Then, for any
given positive integer $e$, there exists an $r$-tuple of
convergent power series ${\cal T}^e(x,u)\in
(\C\{x,u\})^r$ such that $\varphi ({\cal
T}^e(x,u);x,u,y)=0$ in
$\C[[x,u,y]]$ and such that ${\cal T}^e(x,u)$
agrees up to order $e$ (at 0) with ${\cal T}(x,u)$. 
\end{lemma} 
 \begin{proof} First observe
that ${\cal T}(x,u)$ is a formal power series solution of the {\it analytic} system in the unknown $W$,
\begin{equation}\label{eqsyst}
\displaystyle
\left[\frac{\partial^{|\beta|}\varphi}{\partial y^{\beta}}
(W;x,u,y)\right]_{y=x}\equiv 0,\quad \beta \in \N^{q}.
\end{equation} Thus, an
application of Artin's approximation theorem \cite{A} gives, for any
positive integer $e$, an $r$-tuple of convergent power series
${\cal T}^e(x,u)\in (\C\{x,u\})^{r}$ solution in $W$ of
(\ref{eqsyst}), and which agrees up
to order $e$ with ${\cal T}(x,u)$.  The Lemma follows by noticing that $\varphi ({\cal
T}^e(x,u);x,u,y)\equiv 0$ in $\C[[x,u,y]]$ if and only if ${\cal
T}^e(x,u)$ is solution of (\ref{eqsyst}). The
proof of Lemma \ref{lem3} is complete. \end{proof}   The
following proposition will also be useful in the proof
of Theorem
\ref{th1} (see
\cite{BM} for a proof for instance).
\begin{proposition}\label{prop9} Let ${\cal J}(x) =({\cal
J}_1(x),\ldots,{\cal J}_r(x))\in (\C\{x\})^r$,
$x\in \C^k$, $k,r\geq 1$, ${\cal J} (0)=0$, and
${\cal V}(t)\in \C[[t]]$, $t\in
\C^r$.  If ${\cal V}\circ {\cal J}$ is convergent and ${\cal J}$
is generically submersive, then ${\cal V}$ itself is convergent. 
\end{proposition}

\begin{proof}[Proof of Theorem \ref{th1}.] 
Restricting the identity (\ref{eq4000}) to
the set
$$\{(v_3(z',\xi,\eta),\bar{v}_2(\xi,\eta)):(z',\xi,\eta)\in
(\C^{3n-3},0)\} \subseteq \M,$$ where $v_j$, $j=2,3$, are
the Segre sets mappings defined by (\ref{market}), we obtain
\begin{equation}\label{eq7} {\cal R}(v_3(z',\xi,\eta),(\bar
f'\circ
\bar{v}_2)(\xi,\eta))=(\bar{f}_n\circ \bar{v}_2)(\xi,\eta). 
\end{equation}  Here, ${\cal R}$ is the formal power series
defined by (\ref{eqR}). We would like to apply
Lemma
\ref{lem3} to the formal equation  (\ref{eq7}) with $y=z'$,
$x=\eta$, $u=\xi$, ${\cal T}(x,u)=(\bar{f}\circ
\bar{v}_2)(\xi,\eta)$, $W=(\lambda,\mu)$, $\lambda \in
\C^{n-1}$, $\mu \in \C$ and
$$\varphi ((\lambda,\mu);\eta,\xi,z'):={\cal
R}(v_{3}(z',\xi,\eta),\lambda)-\mu.$$  For this, one has to
check that any derivative of the formal holomorphic power
series
$\C^{4n-4}
\ni (z',\xi,\eta,\lambda)\mapsto {\cal
R}(v_3(z',\xi,\eta),\lambda)$ with respect to $z'$ evaluated
at $z'=\eta$ is convergent with respect to the variables
$\lambda$, $\xi$ and
$\eta$.  All these derivatives involve derivatives of $v_3$
at $z'=\eta$ (which are convergent) and derivatives of the form
$\left[{\cal
R}_{z^{\gamma}}(v_3(z',\xi,\eta),\lambda)\right]_{z'=\eta}
$, for $\gamma
\in \N^n$.  Because of the reality condition (\ref{eqm1})
and the definition of $v_1$ and $v_3$ given by
(\ref{market}), we have
$$v_3(\eta,\xi,\eta)=v_1(\eta).$$ This implies that for
each $\gamma \in \N^n$, we have $$\left[{\cal
R}_{z^{\gamma}}(v_3(z',\xi,\eta),\lambda)\right]_{z'=\eta}={\cal
R}_{z^{\gamma}}(v_1(\eta),\lambda),$$ with the
right-hand side being convergent in
$(\eta,\lambda)$ by Proposition
\ref{prop1}.  Thus, by Lemma \ref{lem3}, we have,
for any positive integer
$e$, a convergent power series mapping, denoted
${\cal T}^e(\xi,\eta)=({\cal T}^{'e}(\xi,\eta),{\cal
T}_n^e(\xi,\eta))\in \C^{n-1}\times \C$, which agree up to order
$e$ with $(\bar{f}\circ \bar{v}_2)(\xi,\eta)$ and such that
\begin{equation}\label{eq8} {\cal R}(v_3(z',\xi,\eta),{\cal
T}^{'e}(\xi,\eta))={\cal T}_n^e(\xi,\eta),\ {\rm in}\ \C[[z,\xi,\eta]].
\end{equation} Since ${\cal T}^e(\xi,\eta)$ is a convergent power
series mapping, in order to show that ${\cal R}(z,\lambda)\in
\C\{z,\lambda\}$, it suffices to show by Proposition \ref{prop9}
that for $e$ large enough the generic rank of the holomorphic map 
\begin{equation}\label{eqf}
(\C^{3n-3},0)\ni (z',\xi,\eta)\mapsto
(v_3(z',\xi,\eta),{\cal T}^{'e}(\xi,\eta))\in \C^{2n-1}
\end{equation}
is
$2n-1$. For this, note  that since $M$ is minimal at 0, the map
$\bar{v}_2$ is of generic rank $n$, and thus, by the form of $v_3$
given in (\ref{market}), the holomorphic map
$(\C^{3n-3},0)\ni (z',\xi,\eta)\mapsto
(v_3(z',\xi,\eta),\bar{v}_2(\xi,\eta))\in {\cal M}$ 
is of generic rank $2n-1$
(see \cite{BERb}). Moreover, since $f$ is invertible, we have
${\rm det}\left( \frac{\partial f'}{\partial
z'}(0)\right)\not=0$, which implies that the rank of the formal
map ${\cal M}\ni(z,w)\mapsto  (z,\bar f'(w))$ is $2n-1$ (at the
origin). From this, we see that the formal map
$(\C^{3n-3},0)\ni (z',\xi,\eta)\mapsto (v_3
(z',\xi,\eta),(\bar{f}'\circ
\bar{v}_2)(\xi,\eta))$ has rank $2n-1$. (The rank of
such a formal map is its rank in the quotient field
of $\C[[z',\xi,\eta]]$.)  Since ${\cal
T}^e(\xi,\eta)$ agrees up to order $e$ with
$(\bar{f}\circ \bar{v}_2)(\xi,\eta)$, we obtain that for $e$ large
enough, the mapping (\ref{eqf}) is of generic rank $2n-1$.  This
completes the proof of Theorem \ref{th1}.\end{proof} 
\begin{remark}  When the target hypersurface $M'$ is given in normal
coordinates i.e. $\bar\Phi'(\omega,0)\equiv \omega_n$, then the
normal component $f_n$ of a formal biholomorphism
$f:(M,0)\rightarrow (M',0)$ is convergent provided that the source
hypersurface $M$ is minimal. Indeed, this follows by taking
$\lambda=0$ in Theorem
\ref{th1}. 
\end{remark}
\begin{proof}[Proof of Theorem $\ref{cor2}$.] 
By the Taylor expansion (\ref{eqts}) and by Theorem
\ref{th1}, we obtain that all the
$\phi'_{\alpha}\circ f$ are convergent in a common
neighborhood $U$ of $0\in \C^n$. Since
$M'$ is holomorphically nondegenerate, by \cite{BR, St}, there
exists
$\phi'_{\beta^1}(\omega),\ldots,\phi'_{\beta^n}(\omega)$,
$\beta^i\in \N^{n-1}$, $i=1,\ldots,n$, such that
$${\rm det}\left[\frac{\partial
\phi'_{\beta^i}}{\partial
\omega_j}(\omega)\right]_{1\leq i,j\leq n}\not \equiv
0.$$ Since
$f$ is a formal biholomorphism, this implies that
\begin{equation}\label{eq3001} {\rm det}\left[\frac{\partial
\phi'_{\beta^i}}{\partial
\omega_j}(f(z))\right]_{1\leq i,j\leq n}\not \equiv
0,
\end{equation} as a formal power series in $z$. 
Put
$\psi_{i}(z):=(\phi'_{\beta^i}\circ f)(z)$ and
$R_i(z,\omega):=\phi'_{\beta^i}(\omega)-\psi_i(z)$,
$i=1,\ldots,n$.  Observe that since $\psi_i(z)$ is convergent,
$R_i(z,\omega)\in \C\{z,\omega\}$ for
$i=1,\ldots,n$.  Moreover, since $R_i(z,f(z))=0$, $i=1,\ldots,n$, in
$\C[[z]]$, by (\ref{eq3001}), we may apply
Proposition \ref{prop22} to conclude that $f$ is
convergent.\end{proof}

\section{Transcendence degree and partial convergence of
formal  maps}\label{secpc}
\setcounter{equation}{0}  In this last section, we want to
indicate how Theorem \ref{th1} can be viewed as a result of
partial convergence for formal biholomorphic mappings
of real analytic hypersurfaces. Before explaining
what we mean by this, we need to recall the
following. If $M$ is a
real analytic hypersurface in $\C^n$ and $p\in M$,
let $\K(p)$ be the quotient field of $\C\{z-p\}$,
and
$H(M,p)$ be the vector space over $\K(p)$
consisting of the germs at $p$ of (1,0) vector
fields, with meromorphic coefficients, tangent to
$M$ (near $p$). We then define the {\it degeneracy}
of
$M$ at
$p$, denoted $d(M,p)$, to be the dimension of
$H(M,p)$ over $\K(p)$. It is shown in \cite{BR}
that the mapping $M\ni p \mapsto d(M,p)\in
\{0,\ldots,n\}$ is constant on any connected
component of $M$. Thus, if $M$ is a connected real
analytic hypersurface, one can define its
degeneracy $d(M)$ to be the degeneracy $d(M,q)$ at
any point $q\in M$.  Observe that the germ 
$(M,p)$, $p\in M$, is holomorphically nondegenerate
if and only if
$d(M)=d(M,p)=0$.

Theorem \ref{th1} gives the following result of
partial convergence for formal biholomorphic
mappings of real analytic hypersurfaces. By this,
we mean that we have the following.
\begin{theorem}\label{pc0}
 Let $f : (M,0) \rightarrow (M',0)$ be a formal
biholomorphism between two germs at 0 of smooth
real real-analytic hypersurfaces in
$\C^n$. Assume that $M$ is minimal at $0$ and let
$d(M')$ be the degeneracy of the germ $(M',0)$.
Then, there exists
$g
(\omega)=(g_1(\omega),\ldots,g_{n-d(M')}(\omega))\in
(\C\{\omega\})^{n-d(M')}$, $\omega \in \C^n$, of
generic rank
$n-d(M')$ such that the mapping $g \circ f$ is
convergent.
\end{theorem}
\begin{proof}
 We again use the notations of \S
\ref{sec1}. As in the proof of Theorem \ref{cor2}, we have, using the
expansion (\ref{eqts}),
$$\bar{\Phi}'(f(z),\lambda)=\sum_{\beta
\in \N^{n-1}}(\phi'_{\beta}\circ f)(z)\lambda^{\beta}.$$
Thus,  we know, by Theorem \ref{th1}, that for any multi-index
$\beta \in \N^{n-1}$, $(\phi'_{\beta}\circ
f)(z)$ is convergent in some neighborhood
$U$ of $0$ in $\C^n$.  We choose
$\phi'_{\beta^1}(\omega),\ldots,\phi'_{\beta^r}(\omega)$,
$r=n-d(M')$, of generic maximal rank equal to $n-d(M')$ in a
neighborhood $U'$ of $0$ in $\C^n$ (see \cite{BR}).
Then, if we define
$g_{j}(\omega)=\phi'_{\beta^j}(\omega)$,
$j=1,\ldots,n-d(M')$, we obtain the desired
statement of the Theorem.
\end{proof}
\begin{remark} One should observe that the
convergent power series mapping $g$ in
Theorem \ref{pc0} is obtained in a constructive way
from the target manifold $M'$. Indeed, this is a
consequence of the statement of Theorem
\ref{th1}.
\end{remark}

Now, we want to make explicit links with the notion
of {\it transcendence degree} introduced in
\cite{CPS2} in the
${\cal C}^{\infty}$ mapping problem. For this, we
first set the corresponding definitions in the
formal case.

\begin{definition}\label{def1} Let ${\cal
H}:(\C^N,0)\rightarrow (\C^{N'},0)$ be a formal (holomorphic)
mapping, and 
$V$ be a complex analytic set through the origin in $\C^N\times
\C^{N'}$. Assume that $V$ is given near the origin in $\C^{N+N'}$
by $$V=\{(x,y)\in \C^N\times
\C^{N'}:b_1(x,y)=\ldots=b_q(x,y)=0\},$$
$b_i(x,y)\in \C\{x,y\}$, $i=1,\ldots,q$. Then, the graph of
${\cal H}$ is said to be {\it formally contained} in $V$ if $b_1(x,{\cal H}(x))=\ldots=b_q(x,{\cal H}(x))=0$ 
in $\C[[x]]$. \end{definition}
It follows from the
Nullstellensatz that this definition is independent of the choice
of the defining functions $(b_i)$ for $V$. 

\begin{definition}\label{def2}  Let ${\cal
H}:(\C^N,0)\rightarrow (\C^{N'},0)$ be a formal holomorphic
mapping. Let $V_{\cal H}$ be the germ of the complex
analytic set through the origin in $\C^{N+N'}$ defined as the
intersection of all the complex analytic sets through
the origin in $\C^{N+N'}$ which formally contain the graph of
${\cal H}$. Then the {\it transcendence degree} of ${\cal H}$ is
the nonnegative integer ${\rm dim}_{\C}\ V_{\cal H}-N$.
\end{definition} This definition is motivated by the following
result.
\begin{proposition}\label{prop51} Let ${\cal
H}:(\C^N,0)\rightarrow (\C^{N'},0)$ be a formal holomorphic
mapping. Then, the following conditions are
equivalent:\\ i) ${\cal H}$ is convergent.\\ ii) The
transcendence degree of ${\cal H}$ is zero.
\end{proposition}
\begin{proof} The implication i)
$\Rightarrow$ ii) is clear. The other implication is
equivalent to the following proposition.\end{proof}
\begin{proposition}\label{propn} Let ${\cal
H}:(\C^N,0)\rightarrow (\C^{N'},0)$ be a formal (holomorphic)
mapping. If there exists a germ at $0$ of a
complex analytic set
$V\subset
\C^{N+N'}$ which formally contains the graph of ${\cal H}$
with
${\rm dim }_{\C}V=N$, then ${\cal H}$ is convergent.
\end{proposition} \begin{proof}  Let $V$ be
as in the proposition.  We can assume that, near
$0\in \C^{N+N'}$, $$V=\{(x,y)\in \C^N\times
\C^{N'}:b_1(x,y)=\ldots=b_{p}(x,y)=0\},$$ where each $b_j(x,y)\in \C\{x,y\}$.  To this complex analytic set
$V$, as is customary, we associate the following ideal of
$\C\{x,y\}$ defined by $${\cal I}(V)=\{s\in \C\{x,y\}:s\
{\rm vanishes\ on}\ V\}.$$ By the Noetherian property, we
can assume that ${\cal I}(V)$ is generated by a family
$(h_i (x,y))_{i=1,\ldots,k}\subseteq \C\{x,y\}$.  Furthermore, by the Nullstellensatz
\cite{L, T}, for $i=1,\ldots,k$, there exists an integer $\mu_i$
such that $h_i^{\mu_i}(x,y)$ is in the ideal
generated by the $b_j(x,y)$, $j=1,\ldots,p$.  This
implies that for $i=1,\ldots,k$,
$h_i(x,{\cal H}(x))=0$ in $\C[[x]]$.  The following result, a
consequence of the Artin approximation theorem, is contained in
\cite{T} (p.63).   Assume that the height of the ideal ${\cal I}(V)$
(generated by the family
$(h_j(x,y))_{1\leq j\leq k}$ in $\C\{x,y\}$) is equal to
$N'$. Then, any formal solution ${\cal Y}(x)\in (\C[[x]])^{N'}$, ${\cal
Y}(0)=0$, of the system 
$h_1(x,y)=\ldots=h_k(x,y)=0$ (in the unknown $y$) is convergent. 
Thus, to obtain the convergence of our original formal power series
${\cal H}$, it suffices to check that the height
of ${\cal I}(V)$ is equal to $N'$.  Since
$\C\{x,y\}$ is a local regular ring of Krull dimension $N+N'$, by
Proposition 6.12, p.22 of \cite{T}, we have the formula
$${\rm height}\, ({\cal I}(V))+{\rm dim}\ \C\{x,y\}/{\cal
I}(V)=N+N',$$ where
${\rm dim\ }\C\{x,y\}/{\cal I}(V)$ is the Krull dimension of
the ring
$\C\{x,y\}/{\cal I}(V)$.  Since the Krull dimension of
such a ring coincides with the dimension of the complex analytic
set ${V}$ ({\it cf.}  \cite{L}, p.226, Proposition 1), which is,
here, equal to $N$, we obtain that the height of ${\cal I}(V)$
is $N'$.  This completes the proof of Proposition
\ref{propn}, and hence, the proof of Proposition \ref{prop51}.
\end{proof}

With these tools at our disposal, we can now state
a result which follows from Theorem \ref{pc0}.
\begin{corollaire}\label{thfin} Let $f : (M,0)
\rightarrow (M',0)$ be a formal biholomorphism
between two germs at 0 of smooth real real-analytic hypersurfaces in
$\C^n$. Assume that $M$ is minimal at $0$ and denote by ${\cal
D}_f$ the transcendence degree of the map $f$. Then, ${\cal
D}_f\leq d(M')$, where $d(M')$ is the degeneracy of $M'$. In other
words, there exists a complex analytic set of (pure) dimension
$n+d(M')$ which formally contains the graph of $f$. 
\end{corollaire}
\begin{proof} By Theorem \ref{pc0}, there exists a
convergent power series mapping $g
(\omega)=(g_1(\omega),\ldots,g_{n-d(M')}(\omega))\in
(\C\{\omega\})^{n-d(M')}$ such that for each
$j=1,\ldots,r$, $\delta_j(z):=(g_j
\circ f)(z)$ is convergent, $r=n-d(M')$. Then, the
graph of
$f$ is formally contained in the complex analytic set
$$A=\{(z,\omega)\in
(\C^{2n},0):g_{1}(\omega)-\delta_{1}(z)=\ldots
=g_{r}(\omega)-\delta_{r}(z)=0\}.$$ Let
$A=\cup_{i=1}^k\Gamma_i$ be the decomposition of
$A$ into irreducible components.  For any positive integer
$\sigma$, one can find, according to the Artin approximation
theorem \cite{A}, a convergent power series mapping
$f^{\sigma}(z)\in (\C\{z\})^n$ defined in some small
 neighborhood $U^{\sigma}$ of $0$ in $\C^n$, which
agrees with
$f(z)$ up to order
$\sigma$ (at 0) and such that the graph of $f^{\sigma}$, denoted
$G(f^{\sigma})$, is contained in $A$.  Since $G(f^{\sigma})$ is
contained in $A$, it must be contained in an irreducible
component of $A$.  Thus, by the pigeonhole principle, at least
one subsequence of
$(f^{\sigma})_{\sigma \in \N^*}$ is contained in one of such
irreducible components, say $\Gamma_1$.  There is no loss of
generality in assuming that such a subsequence is
$(f^{\sigma})_{\sigma \in \N^*}$ itself.  We first observe that
this implies that the graph of $f$ is formally contained in
$\Gamma_1$.  Moreover, since
$f$ is a formal biholomorphism, the family $(f^{\sigma})_{\sigma
\in \N^*}$ is also a family of local biholomorphisms.  In
particular, this implies that the generic rank of the family of
holomorphic functions
$$((g_{i}\circ f^1)(z))_{1\leq i\leq r},$$
is
$r$.  As a consequence, if $z_0$ is close enough to
0 in $\C^n$ and is chosen so that the rank of the
preceding family at $z_0$ equals
$r$, the implicit function theorem shows that $A$ is an
$n+d(M')$-dimensional complex submanifold near
$(z_0,f^{1}(z_0))\in
\Gamma_1$.  Since $\Gamma_1$ is irreducible, it is pure-dimensional;
thus
$\Gamma_1$ is an $n+d(M')$ pure-dimensional complex analytic set
formally containing the graph of $f$.  By definition of the
transcendence degree, this implies that ${\cal D}_f\leq
d(M')$.\end{proof}

The following example illustrates the applications
of Theorem \ref{pc0} and Corollary \ref{thfin}.
\begin{example}
Let $M=M'$ be the minimal real algebraic
hypersurface through the origin in $\C^3$ given by 
$${\rm Im}\ z_3=|z_1z_2|^2.$$
Here, $M$ is holomorphically degenerate and its
degeneracy $d(M)$ is equal to 1. Consider the
following formal biholomorphic self-map of $M$:
$$f_h:\C^3\ni (z_1,z_2,z_3)\mapsto
(z_1e^{h(z)},z_2e^{-h(z)},z_3)\in \C^3,$$
where $h(z)=h(z_1,z_2,z_3)$ is any non-convergent
formal power series vanishing at the origin. Observe
that Theorem \ref{th1} gives in this example that
for any formal biholomorphic self-map
$f=(f_1,f_2,f_3)$ of
$M$, the product $f_1f_2$ and the third component
$f_3$ are necessarily convergent. Observe in this example that
the first two components of
$f_h$ are not convergent, but that the transcendence degree
of the map $f_h$ is actually $1=d(M)$. Indeed, the graph of $f_h$
is formally contained in the complex analytic set of
dimension 4
$$V=\{(z,\omega)\in \C^3\times
\C^3:\omega_1 \omega_2=z_1z_2,\
\omega_3=z_3\},$$
and  cannot be formally contained in a complex analytic set
of dimension 3, since otherwise 
$f_h$ would be convergent by Proposition \ref{propn}.
\end{example}

 We
conclude by observing that Theorem \ref{cor2} can
be regarded as a  direct consequence of Corollary
\ref{thfin}. Indeed, when $M'$ is holomorphically
nondegenerate, as mentioned above,
$d(M')=0$ and hence ${\cal D}_f=0$ by Corollary
\ref{thfin}. It then follows from Proposition
\ref{prop51} that $f$ is convergent in that case. 
We should also mention that Theorem \ref{cor2},
Theorem \ref{pc0} and Corollary \ref{thfin} are all
consequences of our main result, namely Theorem
\ref{th1}.
\section*{Acknowledgments}  I would like to thank Makhlouf Derridj for his
interest in this work and Vincent Thilliez for calling my
attention to the paper \cite{N}.  I wish also to address special
thanks to Salah Baouendi and Linda Rothschild for many
simplifying remarks and suggestions.

\end{document}